---

# INNOVATION OF PARAMETRIC PLANE CURVE AND ALGEBRAIC SURFACES FOR BATIK's MOTIFS

H.A Parhusip
Math. Departement
Science and Mathematics Faculty
Universitas Kristen Satya Wacana
Salatiga
hanna.arini@staff.uksw.edu

**Abstract**

Innovation on batik's motifs are shown here. Mathematical approach such as parametric plane curves and algebraic surfaces are introduced to be the novelty of this research. The specific parametric plane curve used in this paper is the modified hypocycloid curve with the MATLAB program. The quadric surface is the particular algebraic surface employed in the Surfer program and varied into several surfaces. Furthermore, the gender issue on batik's painting is examined by introducing Surfer Batik Workshop for 24 junior students. Though statistically the used number does not representative enough, more colors are used by girls than boys, the interest of painting batik's motifs is greater in girls due to the present attitude of participants.

## INTRODUCTION

Batik has been known as an Indonesian heritage with fully art works which mainly denote as motifs textiles in traditional ways of motifs printing or drawn manually. Different city has different type of batik characterizing its culture. Hence various kinds of batik's motifs may be developed from local cultures. Up to now, cultures are mixed unavoidable yielding new type batik's motifs including modern styles in presenting new development in batik, i.e. using abstracts paintings, using different types of natural coloring and involving science and mathematics.

In mathematics particularly, fractal has been included as batik's motifs where various kinds have been readily published in the batik's markets in Indonesia. Hence, batik's motifs generated by fractals have been developed in Indonesia to promote mathematical creations in patterns for ordinary people through textiles. However, improvement for batik's motifs is still growing. As one of these directions, new motifs called BATIMA (*Batik Inovasi Matematika*) are defined based on mathematical curves and surfaces, i.e. most motifs have their own formulas. Unlike fractals, the motifs are created from innovation of parametric curves and or algebraic surfaces with software Surfer developed by Imaginary (research group from Germany) where this software usually is used only for mathematical visualization. Thus originality of this research addresses on developing batik's motifs with that aspects. Since drawing batik is more female's interest, a gender issue arises during the creation of BATIMA. In traditional scale, one has researched this aspect (Dwiyanto and Nugrahani, 2002), empowerment of female batik workers becomes necessary as proposed by Pinta (2013).

An observation due to gender aspect on early ages is done in this research by introducing Surfer batik workshop for junior schools, in 10-16 June 2017. Surfer refers to software developed by Imaginary from Germany which is a research group for mobilizing mathematics. Thus students use Surfer to have batik's motifs and then drawing into textiles with batik materials. The number of students is 24 containing 8 boys and 16 girls. Statistically, this number does not represent the study in the issue of gender interest in drawing batik. However, the research has shown boys have worst attitudes in drawing and finishing the Surfer batik workshop due to less interest than girls. The detail observation is explained below.

**THEORETICAL APPROACH**

**Batik exploration for Gender Equality from junior high school**

Equating males and females access on many areas can be introduced in early ages for children in schools. Students are encouraged to have the same possible activities though natural behavior of each student may affect his/her contribution in activities. For instance, football is much familiar for boys than girls, cooking is more enjoyable for girls than boys naturally.

In schools, varies of activities can be attended for students to explore their capabilities. However natural interest due to gender issue may drive students in the same behaviors. Hence a feminism activity that open with equal access for boys and girls may have girls more participants and vice versa. This issue is addressed here by giving Surfer Batik Workshop activity called Innoweek Junior Summer Camp in 10-16 June 2017 for 24 students from Jakarta and Cirebon. Though the survey has not presented the conclusion above statistically, girls are much interest in drawing batik than boys. Fig. 1 shows one of the documents in this activity.

Not only the interest for attending feminism activity leading to the number of boys less than girls, the esthetics of this activity presents different values due to gender issue. The paper will discuss the results in the section of Result and Discussion.

However, various subjects in mathematics contain patterns are not yet implemented in textiles. Therefore, the research here provides novelty in using calculus and algebraic surfaces as batik motifs where are normally used only in regular teaching and learning mathematics. Visualizations of the motifs from calculus such as parametric curves particularly have been studied in previous researches (Parhusip,2014)(Purwoto,et.al.,2014)(Suryaningsih,et.al.,2013).

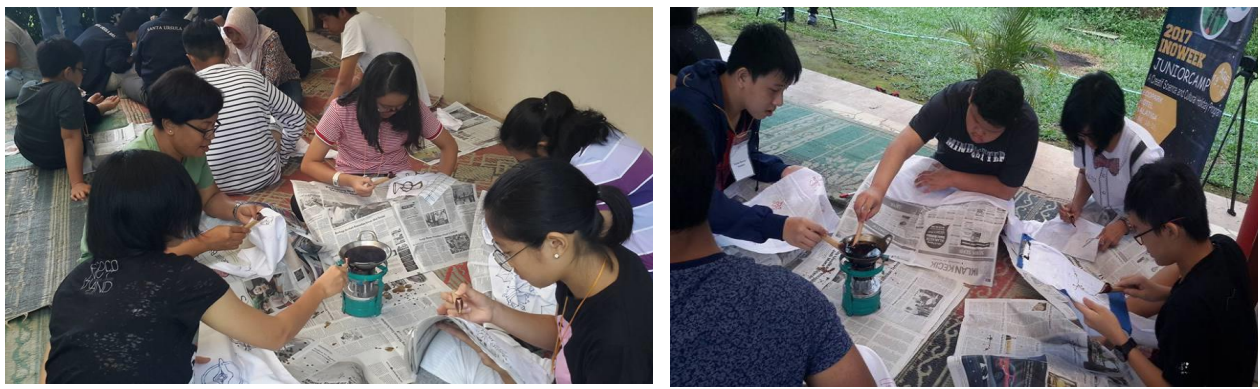

**Figure 1**. Students from Santa Ursula BSD Tangerang and SMP BPK Penabur Cirebon drawing batik in Innoweek Junior Camp, 10-16 June 2017 , Agrowisata Salib Putih Ecopark, Salatiga, Central Java .

The purposes of these researches are creating several motifs in 2 dimensional and 3 dimensional where some patterns show some motifs are similar to natural geometries. The patterns are then printed into bags and t-shirts, and some patterns are considered to be mathematical ornaments (Parhusip,2015). Additionally, the known algebraic surfaces are modified such that new kinds motifs appear, i.e. 3 dimensional motifs in batiks motifs by using modified parametric curve into 2 dimensional curves and 3 dimensional surfaces (Parhusip,2016). Additional software called Surfer is also used where several algebraic surfaces are explored.

**Parametric Curves for Batik Motifs**

There are many parametric curves introduced in calculus for learning calculus in undergraduate mathematics, sciences and engineers. Very few Indonesian have established calculus (particularly) in industry. In developing country such as Indonesia, interaction between mathematics and textile industries is not as high as in developed countries. Therefore, one has industrial mathematics in Indonesia in the lever of home industries.  This is done by using simple mathematics such calculus, parametric curves particularly. As we have known, a parametric curve defined by 2 equations where each point $(x(\theta), y(\theta))$ defined the curve. Frequently, parametric curve contains more than 1 parameters. In the previous research, each point defining the modified hypocycloid curve (Parhusip,2015) contains 2 parameters, i.e :

$$x = x(\theta, a, b), y = y(\theta, a, b) . \qquad (1)$$

Since values of parameters are varied, the curves behave differently yielding various curves. Furthermore, the exploration of this curve through a spherical coordinate leads into a surface, i.e.

$$x = x(\theta,a,b)\rho\sin\phi, \; y = y(\theta,a,b)\rho\sin\varphi, \; z = \rho\cos\phi \quad (2)$$

with

$$\rho = r/\sin\phi, \quad r = \sqrt{x(\theta,a,b)^2 + y(\theta,a,b)^2}.$$

By squaring of each equation and adding the results, one has

$$x^2 + y^2 + z^2 = x^2(\theta,a,b)\rho^2\sin^2\phi + y^2(\theta,a,b)\rho^2\sin^2\varphi + \rho^2\cos^2\phi$$
$$= \rho^2\sin^2\phi\left(x^2(\theta,a,b) + y^2(\theta,a,b)\right) + \rho^2\cos^2\phi. \quad (3)$$

Since $\rho = r/\sin\phi$, Eq.(3) becomes

$$x^2 + y^2 + z^2 = r^2\left(x^2(\theta,a,b) + y^2(\theta,a,b) + \cot^2\phi\right). \quad (4)$$

For a given $(\theta^*,\phi^*)$ and (a,b), each point(x,y,z) in Eq.(4) satisfies its own sphere equation. Thus for any given parametric curve satisfies Eq.(1) can be extended into 3 dimensional case using Eq.(2). Many surfaces are created by Eq. (2) for designing mathematical ornaments and souvenirs (Parhusip,2015) . Furthermore, additional plane curves are innovated by considering Eq.(1) as the domain of the complex mappings to Eq.(1) (Suryaningsih,et.al,2013).

**Example 1.** Let us consider

$$x = x(\theta,a,b) = -(a+b)\sin\theta - (a+b)\sin\left(\left(\frac{a+b}{b}\right)\theta\right),$$

$$y = (a+b)\cos\theta + (a+b)\cos\left(\left(\frac{a+b}{b}\right)\theta\right). \quad (4)$$

The plane curve and its surface are depicted in Fig.2.

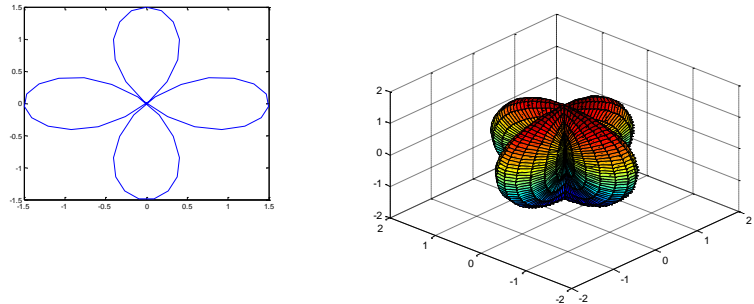

**Figure 2.** Illustration of a plane curve satisfies Eq.(4) (left) and its surface(right).

Both can be one of the innovations for batik's motifs by removing the axes and the grid. Furthermore, mathematics transformations such as rotation, shearing, scaling are available to innovate the original curve or surface to be new curve and surface.

**Algebraic Surfaces for Batik Motifs**

We have known that a surface in 3 dimensional takes in the form *f(x,y,z)*= 0. One of the typical algebraic surfaces is a quadratic surface with its general form is

$$a_1 x^2 + a_2 y^2 + a_3 z^2 + 2a_4 yz + 2a_5 zx + 2a_6 xy + 2a_7 x + 2a_8 y + 2a_9 z + a_{10} = 0.$$

We have some examples in this family, for instance a cone, cylinder, sphere, paraboloid, hyperbolic paraboloid. Quadric surface is denoted by the maximum sum of powers of all terms in $a_m x^{i_m} y^{j_m} z^{k_m}$ is 2. Since the maximum of *i+j+k* may be > 2, one has higher degree of surfaces classifying algebraic surfaces into cubic, quartic, quintic, sextic and septic surfaces for the degree of surfaces are 3,4,5,6,7 respectively. Higher degrees are still possible.

Cubic surface contains a family of surfaces with the maximum of *i+j+k* =3 in the polynomial equation, i.e.

$$a_1 x^3 + a_2 y^3 + a_3 z^3 + a_4 x^2 y + a_5 x^2 z + a_6 xy^2 + a_7 xz^2 + a_8 y^2 z + a_9 z^2 y$$
$$+ a_{10} xyz + a_{11} xy + a_{12} yz + a_{13} x + a_{14} y + a_{15} z + a_{16} = 0.$$

The equations sometimes are written parametrically. However drawing the algebraic surfaces in the used software here (Surfer 2008) is not available for using parametric equation. Several surfaces in higher degree have been collected into world recorded surfaces characterized by many singularities. In this research, students in Surfer Batik Workshop are allowed to draw in their own creativities such that the surfaces may differ variously than the equations above.

**Algebraic Surfaces through Surfer**

Surfer is one of software for creating surfaces easily, i.e. we only type the algebraic equations and Surfer does drawing directly as soon as the common is typed properly. In order to understand how Surfer is working, some initial steps are explained here. One has known, a typical algebraic surface is a sphere, i.e. $x^2 + y^2 + z^2 = a^2$. Surfer requires that the surface is written as *f*=0, hence $x^2 + y^2 + z^2 - a^2 = 0$. Additionally, Surfer allows us to create several surfaces. Therefore, the general formulas in writing an algebraic surface in Surfer must follow $\prod f_i = 0$.

The readers are expected to have knowledge on a circle equation with a center O and its radius $r$ on ($x,y$)-plane , i.e.

$$x^2 + y^2 = r^2. \qquad (5)$$

Similarly, the equation for a circle with its center coordinate is ($h,k$) and its radius $r$

$$(x-h)^2 + (y-k)^2 = r^2. \qquad (6)$$

One may have the similar expressions on ($y,z$) and ($x,z$) planes. If we take the Eq.(5) and Eq.(6) are true for all values $z$ in 3 dimensional space ($x,y,z$), these equations express a cylinder for all values of $z$. For instance, $x^2 + y^2 = 1$ or $x^2 + y^2 - 1 = 0$ describes a cylinder for all values of $z$. Using Surfer, one simply writes the Eq.(5) as x^2+y^2-1=0 shown in Fig.2.

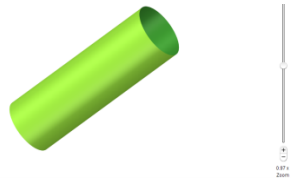

**Figure 2.** Illustration of Eq.(5) using Surfer with vertical line is the menu for zooming.

By adding another variable into Eq.(5), we will have other quadric surface called a sphere, i.e $x^2 + y^2 + z^2 = r^2$. Now we will use the ball equation to create an innovative surface which can be used as batik motif.

**Example 2.** The purpose of the processes here is creating 2 balls and make both into Poke Planet ball (a well known cartoon). Using several times of the sphere equation, one has the Poke Planet Ball with the Equation in Surfer language as

((x^2+y^2+z^2-100))*((x^4+y^4-2)*(x^3+y^3-4)*(x^4+y^4-8)*(x^3+y^3-16)*(x^4+y^4-32)*(x^3+y^3+36)) .

The resulting Poke Planet Ball is drawn in Fig. 3 using the zoom button into 0.04x.

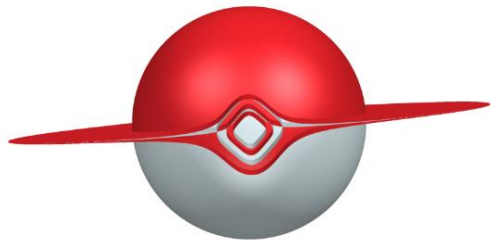

**Figure 3.** Poke Planet Ball surface

**Example 3.** One uses the menu with option Ding dong (*f*=0),i.e. (x^2+y^2+z^3-z^2)=0 (shown by Fig. 4.left) and then create another surface (*g*=0), i.e. (((x-y)*(x+y)-a)*((x+y)*(x-y)+a))=0 (depicted in Fig.4 middle). Combining both equations means that both surfaces into one surface by setting *f g* = 0. The Surfer requires to write it into

(x^2+y^2+z^3-z^2)*(((x-y)*(x+y)-a)*((x+y)*(x-y)+a))=0.

The parameter *a* will be a new menu automatically coming to the system such that *a* is freely chosen. Employing the zoom menu 0.91x and *a*=0.02, an innovative surface is built called the Atom fish.

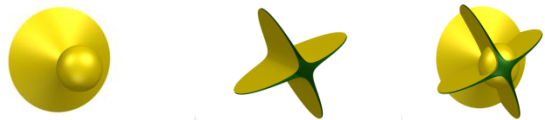

**Figure 4.** Atom Fish surface (right)

**Example 4. Ring Blaster**

As in Example 3, one has to define (*f*=0)in this case, Surfer 2008 is written as (x^2+y^2-1)^2+(y^2+z^2-1)^2-a=0 .

Using the zoom menu into 0.47x and *a*=0,02, the first surface is obtained. The 2nd surface is given by (g=0): (((x-y)*(x+y)-a)*((x+y)*(x-y)+a))=0. Employed again the zoom menu 0.47 and *a*=0.02 ,and hence one has the new innovative surface called Ring Blaster

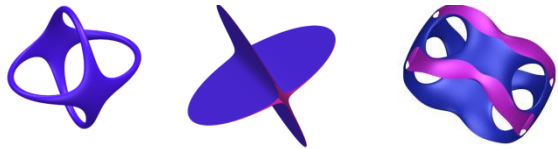

**Figure 5.** Surfaces to create Ring Blaster (left and middle) ; the Ring Blaster (right)

Practically, the obtained surfaces can be explored into several purposes, such as mathematical ornaments, motifs in t-shirts, or motifs in bags or souvenirs designs.

**Drawing Curves as Intersection of Surfaces**

A curve can be considered as intersection of surfaces. As we have known, if we construct 2 surfaces into 1 surface on Surfer 2008, then we write *f g*=0. Let us consider both surfaces intersection yielding a curve, the Surfer 2008 requires us to write the command into $f^2 + g^2 = 0$. As a result, the statements *f*=0 and *g*=0 are satisfied. One can draw the isolated points by adding perturbation parameter in the equation, i.e. $f^2 + g^2 = a$.

**Exxample 5.** Drawing a cross section of 2 tubes

We consider 2 cylinders, i.e. $x^2 + y^2 - 1 = 0$ and $z^2 + (y+3-6b)^2 - 1 = 0$. These cylinders will have cross sections. The cross section gives a curve drawn in Fig. 6.

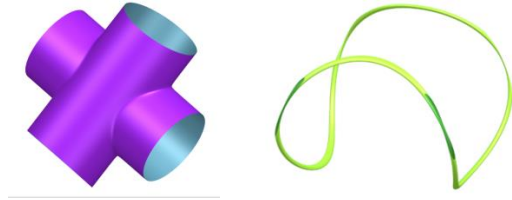

**Figure 6.** Visualization of 2 cylinders (left) drawn by Surfer 2008 and intersection curve (right)
*(x^2+y^2-1)^2+(z^2+(y+3-6\*b)^2-1)^2-0.01\*a=0*
*where a=0.26 and b=0.56, zoom 0.71+.*

However drawing curves into batik's motifs will not be suitable since one has to put the boundary called malam as a traditional material for initial proccess of creating batik's textiles. Therefore the innovation on curves will not be done in this research.

**RESULT AND DISCUSSION**

**Batik from parametric curve**

Several designs have been created due to Eq.(2)-(3). Additional mathematical terms can also be applied to the obtained surface. One example is done by repeating the surface following the Fibonacci sequence. Hence several surfaces can be included into 1 motif as depicted in Fig. 7 and its textile is also shown.

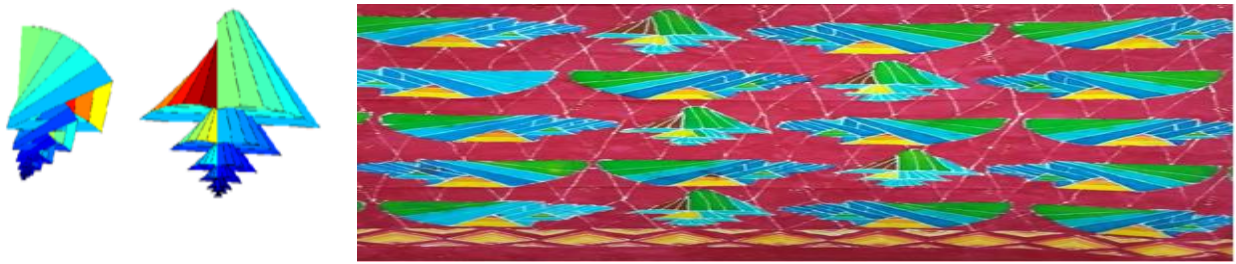

**Figure 7.** Batik's motif (left and middle) obtained from Eq.(2)-(3) and its batik drawn with MATLAB.The textiles of batik is shown (right).

---

**Batik Motif from Algebraic Surfaces**

Each surface on Chapter 2 can be a motif for Batik and hence possible other designs may be created. By employing the same procedure as above, one obtained several motifs. Most of the designs here are based on quadric surfaces.

The surface is created by using $fg = 0$ where $f$ and $g$ are

$$x^2 + y^2 + z^2 + 2xyz - 1 = 0 \text{ and } \left((x-1)^2 + (y-1)^2 + (z-1)^2 + 2(x-1)(y-1)(z-1) - 2\right) = 0$$

respectively. The standard process of manual batik painting is then used into textile. The similar surface is then obtained by replacing $g(x/2, y/2, z/2) = 0$, i.e.

$$\left((x/2)^2 + (y/2)^2 + (z/2)^2 + 2(x/2)(y/2)(z/2) - 1\right) = 0. \tag{7}$$

Another result is done by creating $\prod_{i=1}^{3} g_i = 0$ where $g_1 = 0$ is governed by Eq.(7) and the others are $g(x/k, y/k, z/k)=0$, $k=3,5$. Thus finally motifs and its textile have been made. The motifs and the related batik are depicted in Fig. 8.

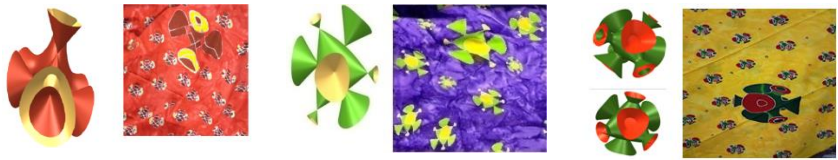

**Figure 8. (left):** Batik's motif (left) and its textile;The motif is obtained by using *fg*=0, i.e.(x^2+y^2+z^2+2*x*y*z-1)*((x-1)^2+(y-1)^2+(z-1)^2+2*(x-1)*(y-1)*(z-1)-2)=0.
 (**middle**): The motif is obtained by using *fg*=0, (x^2+y^2+z^2+2*x*y*z-1)*((x/2)^2+(y/2)^2+(z/2)^2-2*x/2*y/2*z/2-1)=0.
 (**righ**t):; The motif is obtained by using *fgh*=0,i.e. (x^2+y^2+z^2+2*x*y*z-1)*((x/2)^2+(y/2)^2+(z/2)^2-2*x/2*y/2*z/2-1)*((x/3)^2+(y/3)^2+(z/3)^2-2*x/3*y/3*z/3-1) *((x/5)^2+(y/5)^2+(z/5)^2-2*x/5*y/5*z/5-1)= 0.

**The result of Batik with Surfer's motifs due to gender issue**

In this section, the results of Surfer Batik workshop for boys and girls in 10-16 June 2017 are shown. Fig.9 depicts batik's motifs created by boys and Fig.10 batik's motifs created by girls where each person may draw more than 1.

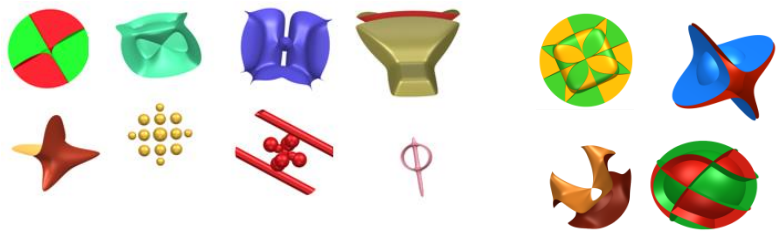

**Figure 9.** Batik's motifs created by boys in Surfer Batik Workshop Innoweek Junior Camp,10-16 June 2017.

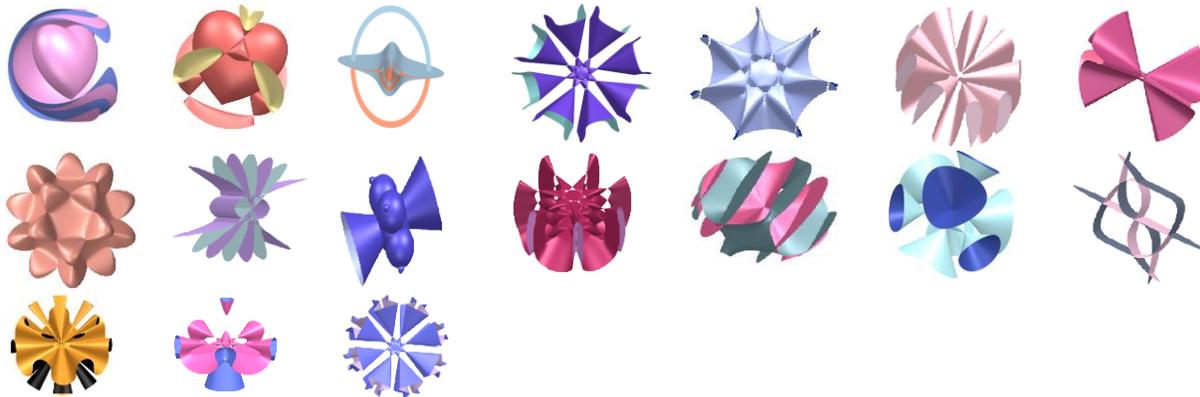

**Figure 10.** Batik's motifs created by girls in Surfer Batik Workshop Innoweek Junior Camp,10-16 June 2017.

There exists no objective measurement to measure the interest of boys less than girls due to the results though some profiles can be concluded to be females like than males by comparing Fig.10 and Fig.11. However, the attitudes during drawing batik's motif may reflects this issue. Additionally, during painting to produce batik, girls showed greater passion than boys. Several related motifs and the textiles are illustrated in Fig.12. Moreover, selected colors due to gender are not observed into detail in this study, though a research on this case has been done by other researches (Black and Wright,2013) where 216 participants were examined and females used more colors than males. One knows that surfer can only provide at most 2 colors in 1 image. The variation of colors are employed in batik's paintings. In the Surfer Batik Workshop, girls applied more colors though only very few students. Thus, the research on this paper has agreed the same conclusion for a gender issue for colors selections in drawing.

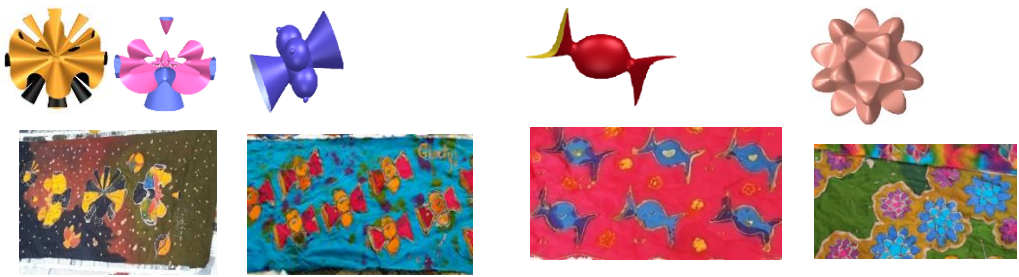

**Figure 12**. Motifs from Surfer (above) and the related batik (below)

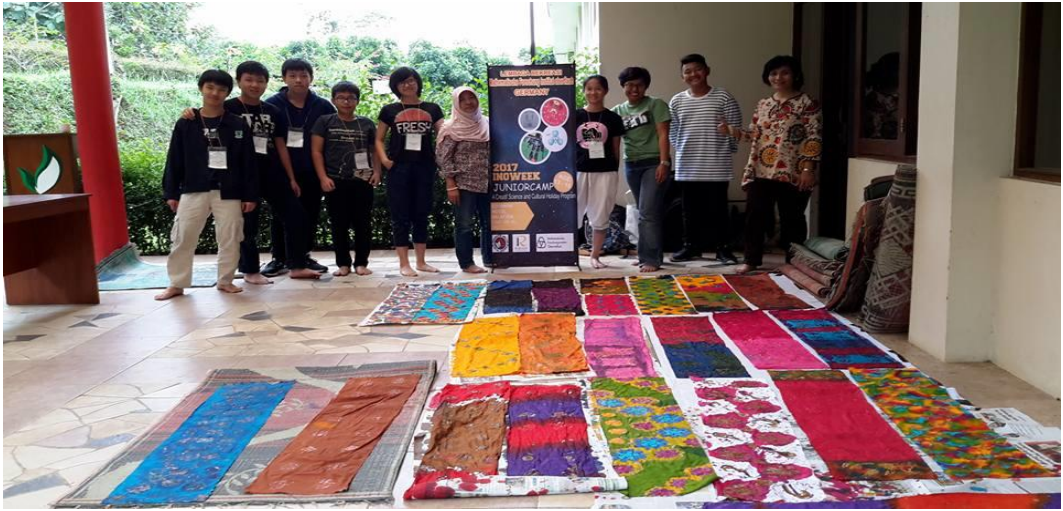

**Figure 12**. Motifs from Surfer (above) and the related batik (below)

**CONCLUSION**

This paper describes innovation on batik's motifs from calculus and algebraic surfaces. The parametric curves from calculus are extended into surfaces through spherical coordinates. Parameters and the number of points used in drawing are varied yielding various types of surfaces to be candidate of motifs. The selected motifs are then implemented into textiles. Similarly, algebraic surfaces are usually known only in mathematics. This research has been successful to innovate algebraic surfaces into batik's motifs.

The gender issue is examined by introducing batik's drawing in Surfer Batik Workshop for juniors school in 10-16 June 2017. The interest of drawing and painting, the used colors have observed leading to greater variants in batik's motifs for girls than boys.

**Acknowledgement**

The research here is funded by Ristekdikti (Ministry of Research, Technology and Higher Education of the Republic of Indonesia) under the grant Penelitian Produk Terapan, No.025/E3/2017 2nd year,2017/2018. Entitled 'Pembuatan GARISMA (Galeri Riset Dan Inovasi Matematika) Berbasis Pemetaan Fungsi Kompleks Pada Kurva Hypocycloid Guna Menarik Minat Belajar Matematika.'
Thank to mathematics students from FKIP UKSW and students from SMP Santa Ursula BSD and BPK Penabur for creating some surfaces (2016) and batik's motifs respectively in this paper.